\documentclass[leqno,oneside,english]{article}
\usepackage[T1]{fontenc}
\usepackage{amsmath}
\usepackage{amssymb}
\usepackage{color}
\usepackage{dsfont}
\usepackage{pstricks}
\usepackage{enumitem}

\usepackage{latexsym}

\newcommand{\RR}{\mathbb{R}}

\newcommand{\ZZ}{\mathbb{Z}}

\makeatletter
\usepackage[T1]{fontenc}

\usepackage{graphicx}

\makeatletter \textwidth15cm \textheight20truecm

\newcommand{\R}{{\mathbb R}}
\newcommand{\C}{{\mathbb C}}

\newcommand{\caT}{{\mathcal T}}

\def\be{\begin{equation}}
\def\ee{\end{equation}}
\def\beq{\begin{eqnarray}}
\def\eeq{\end{eqnarray}}
\def\beqs{\begin{eqnarray*}}
\def\eeqs{\end{eqnarray*}}
\def\ea{\end{array}}
\def\ea{\end{array}}
\def\bt{\begin{thm}}
\def\et{\end{thm}}
\def\br{\begin{rk}}
\def\er{\end{rk}}
\def\bc{\begin{cor}}
\def\ec{\end{cor}}
\def\bl{\begin{lem}}
\def\el{\end{lem}}

\newtheorem{thm}{Theorem}[section]

\newtheorem{cor}[thm]{Corollary}

\newtheorem{lem}[thm]{Lemma}

\newtheorem{rk}[thm]{Remark}

\newenvironment{proof}[1][Proof]{\textbf{#1.} }{\ \rule{0.5em}{0.5em}}

\newcommand{\Om}{\Omega}

\makeatother

\usepackage{babel}
\makeatother
\author{Serge Nicaise\footnote{Universit\'e de Valenciennes et du Hainaut Cambr\'esis, LAMAV,
FR CNRS 2956, Institut des Sciences et Techniques of Valenciennes,
F-59313 - Valenciennes Cedex 9 France,
Serge.Nicaise@univ-valenciennes.fr}, Hengguang Li\footnote{Department of Mathematics, Wayne State University, Detroit, MI 48202, USA, hli@math.wayne.edu}, and Anna Mazzucato,\footnote{Penn State University, University
Park, PA 16802, USA, alm24@psu.edu}}

\begin{document}

\title{Regularity and a priori error analysis of a Ventcel problem in polyhedral domains}

  \maketitle

\begin{abstract}
We consider the regularity of a mixed boundary value problem for the Laplace operator on a polyhedral domain, where Ventcel boundary conditions are imposed on one face of the polyhedron and Dirichlet boundary conditions are imposed on the complement of that face in the boundary. We establish improved regularity estimates for the trace of the variational solution on the Ventcel face, and use them to derive a decomposition of the solution into a regular and a singular part that belongs to suitable weighted Sobolev spaces. This decomposition, in turn, via interpolation estimates both in the interior as well as on the Ventcel face, allows us to perform an {\em a priori} error analysis for the Finite Element approximation of the solution  on  anisotropic graded meshes. Numerical tests support the theoretical analysis.
\end{abstract}

\bigskip



\bigskip

\section{Introduction}

This article concerns the regularity of solutions to an elliptic boundary-value
problem for the Laplace operator on a polyhedral domain in
$\RR^3$ under so-called {\em Ventcel} or {\em Wentzell boundary conditions}. The
regularity result we establish in weighted Sobolev spaces gives rise, in turn,
to {\em a priori} error estimates for the Finite Element Method (FEM) on a
suitable anisotropic mesh.

We first introduce the Ventcel boundary-value problem.
Let $\Omega$ be a bounded domain of $\R^3$ with  Lipschitz boundary $\Gamma$
and let $\Gamma_V$ be an open subset of $\Gamma$
with positive measure.  We  denote by
$\Gamma_D=\Gamma\setminus \Gamma_V$ the complement of $\Gamma_V$, which we assume has also positive measure.

We consider the following mixed boundary-value problem:
\begin{subequations}
\label{Ventcel}
\begin{align}
 -\Delta u &= f,
&\hbox{ in }\Omega,
\label{Ventcel.a}\\
u &= 0  &\hbox{ on } \Gamma_D,
\label{Ventcel.b}\\
-\Delta_{LB} u+\partial_\nu u &=g, &\hbox{ on } \Gamma_V,
\end{align}
\end{subequations}
where $\Delta$ is the standard (Euclidean) Laplacean in $\RR^3$, $\Delta_{LB} $
is the Laplace-Beltrami operator on $\Gamma$, $\nu$ is the unit outer normal
vector on  $\partial \Omega$,  $\partial_\nu $ means the associated normal
derivative, and $f$ and $g$ are given data.

This problem is a special case of a more general boundary-value problem, where
\eqref{Ventcel.b}
is replaced by:
\[
     - \alpha \Delta_{LB} u +  \partial_\nu u +\beta u= g,
\]
which can be thought of as a generalized Robin-type boundary condition. The
more general problem is well posed only under conditions on the sign of
$\alpha$ and $\beta$. Ventcel boundary conditions arise naturally in many
contexts. In the context of multidimentional diffusion processes, Ventcel
boundary conditions were introduced in the pioneering work  of Ventcel
\cite{Ventcel56,Ventcel59} (see also the work of Feller for one-dimensional
processes \cite{Feller52,Feller57}). They can model heat conduction in materials
for which the boundary can store, but not absorb or transmit heat. They can
also be derived as approximate boundary conditions in asymptotic problems
or artificial boundary conditions in exterior problems (see e.g.
\cite{BnHDV09,BnDHV15,NL10} and references therein), in particular in
fluid-structure interaction problems.

Problem \eqref{Ventcel} is known to have a unique variational solution if $f$
and $g$ are in the appropriate Sobolev space as recalled in Section \ref{sreg}.
We are concerned here with the higher regularity for solutions to this problem
when the data is also regular, in the case that the domain $\Omega$ is a
polyhedral domain in $\RR^3$. It is well known that, due to the presence of
edges and corners at the boundary of $\Omega$, even when $\Gamma_V$ is empty,
elliptic regularity does not hold, and the solution is not smooth even if the
data is smooth. This loss of regularity affects the rate of convergence of the
Finite Element approximation to the solution if uniform meshes are used.

By using weighted Sobolev spaces, where the weights are the distance to the
edges and vertices, respectively,  one can characterize precisely the behavior
of the variational solution near the singular set in terms of singular function
and singular exponents  (Theorem \ref{tregu}). In turn, the decomposition
of the solution into a regular and a singular part, together with
interpolation estimates (Theorem \ref{tinterpolation}), leads to establishing
{\em a priori} error estimates for the Finite Element approximation (Corollary
\ref{cerrorestimate}), where the elements are given on an anisotropic mesh that
exploits the improved regularity of the solution along the edges versus the
corners of the polyhedron. There is a well established literature on this
approach for mixed Dirichlet, Neumann, and even standard Robin boundary
condition (see for example \cite{apel:97c, BNZ107, MR3454358}). There are also several works in the literature concerning the
Ventcel boundary-value problems on singular domains (see in
particular \cite{lemrabet:85,PS13}), and their  implementation of the FEM (see
\cite{Kashiwaba:15} and  references therein).
The novelty of this work consists in extending the approach using weighted
spaces and anisotropic meshes to the Ventcel boundary conditions, which
include tangential differential operators at the boundary of the same order
as the main operator in $\Omega$. As a matter of fact, the main difficulty in
considering such boundary conditions lies in establishing the needed regularity
of the traces on the faces of the polyhedron.
For simplicity, we restrict here to the case where the Ventcel condition is
imposed on only one face of the polyhedron. If the Ventcel condition is imposed
on adjacent faces, one would expect higher regularity to hold for the solution
on these faces, under suitable transmission conditions at the common edges.
however, capturing this behavior entails studying weighted Sobolev spaces for
which the weight is the distance to the boundary and not the distance to the
singular set (as those arising from the analysis of equations with degenerate
coefficients). We reserve to address this problem in future works.

The paper is organized as follows.
In Section \ref{sreg}, we recall the variational formulation for Problem
\eqref{Ventcel}, and prove our main regularity result for the solution in
weighted spaces. In Section \ref{sapprox}, we introduce the anisotropic mesh and
the associated Finite Element discretization of the problem, and derive {\em a
priori} error estimates. Section \ref{sinterror2D} contains some refined
2D interpolation estimates valid on the polyhedral faces, needed for the error
analysis. We close in Section \ref{snumerics} by presenting some numerical
examples to validate the theoretical analysis.

We end this Introduction with some needed notation.

If $\Om$ is a domain of $\R^n$, $n\geq 1$, we employ the standard notation
$H^m(\Omega)$ to denote the Sobolev space that consists of functions whose $i$th
derivatives, for $0\leq i\leq m$,  are  square-integrable.   The
$L^2(\Om)$-inner product (resp. norm)
will be denoted by  $(\cdot,\cdot)_\Om$ (resp. $\|\cdot\|_\Om$).   The usual
norm and semi-norm in  $H^{s}(\Om)$, for $s\geq 0$, are denoted by
$\|\cdot\|_{s,\Om}$ and $|\cdot|_{s,\Om}$, respectively.
The trace operator from $H^1(\Om)$ into $H^{\frac12}(\partial \Om)$ will
be denoted  by $\gamma$.
We also introduce the space:
\[
H^1_{\Gamma_D}(\Om)=\{u\in H^1(\Om): \gamma u=0\hbox{ on } \Gamma_D\},
\]
which is clearly a closed subspace of $H^1(\Om)$. If $v$ is a $d$-dimensional
vector, we will write $v\in H^s(\Omega)^d$, although for ease of notation, we
may write $H^s(\Omega)$ simply for $ H^s(\Omega)^d$. Lastly, we employ the
standard notation $\mathcal{D}'(\Omega))$ to denote the space of distributions
on $\Omega$.

Throughout, the notation $A \lesssim B$ is used for the estimate $A \leq C \ B,$
where $C$ is a generic constant that does  not depend on $A$ and $B$.
The notation $A \sim B$ means that both $A \lesssim B$ and $B \lesssim A$ hold.
We will also employ standard multi-index notation for partial derivatives in
$\RR^d$, i.e., $\partial^\alpha=\partial_{x_1}^{\alpha_1}
\ldots\partial_{x_d}^{\alpha_d}$ where $\alpha=(\alpha_1,\ldots,\alpha_d) \in
\ZZ_+^d$ and $|\alpha|= \alpha_1+\ldots+\alpha_d$.

\bigskip

\noindent{\bf Acknowledgements:}  The second author was partially supported by the US National Science Foundation (NSF) grant DMS-1418853.  The third author was partially supported by NSF grant DMS-1312727.
The visit of the first author to Penn State University and Wayne State University, where part of this work was conducted, was partially supported through NSF grant DMS-1312727 and the Wayne State University Grants Plus Program.

\section{Some regularity results\label{sreg}}

In this section we recall needed facts about the well-posedness of the Ventcel
Problem \eqref{Ventcel}, and establish regularity estimates for its variational
solution in weighted spaces.

The variation formulation of \eqref{Ventcel} is well known (see
\cite{AMbook80,lemrabet:85,Kashiwaba:15}).
We let
\[
V:=\{u\in H^1_{\Gamma_D}(\Om): \gamma u \in H^1_0(\Gamma_V)\},
\]
which is a Hilbert spaces equipped with the natural norm
\[
\|u\|_V^2:=\|u\|_{1,\Om}^2+| \gamma u|_{1, \Gamma_V}, \forall u\in V.
\]

We further introduce the bilinear form
\[
a(u,v)=\int_\Om \nabla u\cdot \nabla v\, dx+\int_{\Gamma_V} \nabla_T(\gamma u)\cdot \nabla_T (\gamma v)
\, d\sigma(x), \quad \forall u,v\in V.
\]
As this bilinear form is continuous and coercive in $V$, by the Lax-Milgram
lemma,
for any $f\in L^2(\Om)$ and $g\in L^2(\Gamma_V)$, there exists a unique solution
$u\in V$ of
\be\label{weakform}
a(u,v)=\int_\Om f  v\, dx+\int_{\Gamma_V} g v\, d\sigma(x), \quad \forall v\in V.
\ee

It was shown in \cite[Thm 3.3]{Kashiwaba:15} that if $\Gamma_D$ is empty and if $\Gamma$ is $C^{1,1}$, then
$u$ belongs to $H^2(\Om)$ and  $\gamma u$ belongs to $H^2(\Gamma_V)$.
This regularity is no longer valid if $\Om$ is a non-convex
polyhedral domain, and the main purpose of this section is to describe the
behavior of the solution near the singular set, which consists of the
edges and vertices of the boundary faces of the polyhedron, and characterize
the regularity of boundary traces of $u$ and its derivatives.

To this end, we will employ anisotropic weighted Sobolev spaces, for which the
weights are (variants of) the distance to the edges and vertices, respectively.
There is a vast literature concerning the use of weighted Sobolev spaces in the
analysis of singular domains (we refer for instance to
\cite{dauge:88,KMR1,MazPlamNdim2, MazyaRossmann, NP} and references therein).
In the context of the analysis of Dirichlet/Neumann boundary
conditions, anisotropic Sobolev spaces were used in
\cite{apel:92,apel:97c,Apeletal:14,BNZ07}.

{\em From now on we assume that $\Om$ is a polyhedral domain of the space
and that $\Gamma_V$
is reduced to one face $F$ of the boundary.}

By a face, we mean an open face on the boundary. Let  $\mathcal S$  and
$\mathcal E$ be  the set of vertices and the set of open edges
of $\Omega$, respectively.

On the polygonal face $F$, we denote its set of vertices by $\mathcal S_F$. Given a
vertex $S\in \mathcal S_F$, we denote by $(r_S, \theta_S)$ the radial distance and
angular component of the local polar coordinate system centered at $S$ on the
plane containing $F$. In addition, we let $\omega_{F, S}$ be the
interior angle on the face $F$ associated with the vertex $S$.

Following \cite{apel:97c}, we consider a
triangulation $\{ \Lambda_\ell\}_{\ell=1}^L$
of the domain $\Omega$ that consists of disjoint tetrahedra
$\Lambda_\ell$. We will refer to each tetrahedron $\Lambda_\ell$ as a {\em
macro element}, to distinguish it from the elements of the mesh utilized in the
analysis of the FEM in Section \ref{sapprox}. The purpose of the macro elements
 is to localize the construction and the regularity estimates near edges and
vertices of $\Omega$. We will also refer to any edge or vertex of an element
$\Lambda_\ell$ as a {\em singular edge} or {\em singular vertex}, if that edge or vertex lies along a true edge or is  a true  vertex of $\Omega$ and the solution is not in $H^2$ near that true edge or vertex.

We will assume that each $\Lambda_\ell$ contains at most one singular
edge and at most one singular vertex. If $\Lambda_\ell$ contains both  a
singular edge and a singular vertex, that vertex belongs to that  edge. We will
also assume that all
$\Lambda_\ell$ are shape regular with diameter of order $O(1)$. In each
macro element $\Lambda_\ell$, we introduce a local Cartesian coordinate
system $x^{(\ell)}=(x^{(\ell)}_1, x^{(\ell)}_2, x^{(\ell)}_3)$, such that
the singular
vertex, if it exists, is at the origin,  and the singular edge, if it exists,
lies along the $x^{(\ell)}_3$-axis. We then define the distance functions to
the set of singular edges and singular vertices, respectively, as follows:
\begin{subequations} \label{eqn.dist}
\beq
r^{(\ell)}(x^{(\ell)})&=&\sqrt{(x_1^{(\ell)})^2+(x_2^{(\ell)})^2},\label{eqn.r}
\eeq
\beq
R^{(\ell)}(x^{(\ell)})&=&\sqrt{(x_1^{(\ell)})^2+(x_2^{(\ell)})^2+(x_3^{(\ell)}
)^2},\label{eqn.r1}
\eeq
and introduce the auxiliary function
\beq
\theta^{(\ell)}(x^{(\ell)})&=&r^{(\ell)}(x^{(\ell)})/R^{(\ell)}(x^{(\ell)}
)\label{eqn.r2}.
\eeq
\end{subequations}
We observe that $r^{(\ell)}$, and $R^{(\ell)}$ extend as continuous functions to
the closure of the macro element $\Lambda_\ell$, while $\theta^{(\ell)}$ extends
as a bounded function.

In what follows, we will omit the sup-index $(\ell)$ in these distance
functions and in $x$, when there is no confusion about the underlying macro
element.
Given a subdomain $\Lambda\subset\Omega$, we define the following weighted
Sobolev space for $k\in\mathbb N$ and $\beta,
\delta\in\mathbb R$:
$$
V^{k}_{\beta,\delta}(\Lambda):=\{v\in \mathcal D'(\Lambda), \ \|v\|_{V^{k}_{\beta,\delta}(\Lambda)}<\infty\},
$$
where
$$
 \|v\|_{V^{k}_{\beta,\delta}(\Lambda)}^2=\sum_{|\alpha|\leq k}\|R^{\beta-k+|\alpha|}\theta^{\delta-k+|\alpha|}\partial^{|\alpha|}v\|^2_{L^2(\Lambda)},
$$
and $R(x)=R^{(\ell)}(x^{(\ell)})$ and $\theta(x)=\theta^{(\ell)}(x^{(\ell)})$,
given in (\ref{eqn.r1}) and in (\ref{eqn.r2}), if $x\in \Lambda_\ell$ is
represented by $x^{(\ell)}$ in local coordinates.

We will also need to define spaces on the faces of $\Omega$. To this end, given
 $G$ a bounded polygonal domain in $\mathbb R^2$, we also define
$$
V^k_\gamma(G):=\{v\in \mathcal D'(G), \ \|v\|_{V^{k}_{\gamma}(G)}<\infty\},
$$
where
$$
 \|v\|_{V^{k}_{\gamma}(G)}^2=\sum_{|\alpha|\leq
k}\|\rho^{\gamma-k+|\alpha|}\partial^{|\alpha|}v\|^2_{L^2(G)}.
$$
Above, $\rho$ is the distance function to the set of vertices of $G$, defined
in a manner similar to $R$ above.

We further classify the initial macro elements $\Lambda_\ell$ into four types as
follows:
\begin{enumerate}[label={\bf Type \arabic*.}, align=left, labelwidth=*]
 \item $\bar\Lambda_\ell\cap (\mathcal
S\cup \mathcal E)=\emptyset$;
  \item $\bar\Lambda_\ell\cap \mathcal S\neq \emptyset$ but
$\bar\Lambda_\ell\cap \mathcal E=\emptyset$;
  \item $\bar\Lambda_\ell\cap
\mathcal E\neq \emptyset$ but $\bar\Lambda_\ell\cap \mathcal S=\emptyset$;
   \item $\bar\Lambda_\ell\cap \mathcal E\neq \emptyset$
and $\bar\Lambda_\ell\cap \mathcal S\neq\emptyset$.
\end{enumerate}

We first start with an improved regularity of $\partial_\nu u$ on $\Gamma_V$.
In what follows, for ease of notation we will let $u_F$ be the trace $\gamma
u$ of  $u$ on the face $F$.
Furthermore,
for a two-dimensional domain $D$, we define
the space  $\tilde{H}^s(D)$, $0<s<1$, as the closure of $C^\infty_c(D)$
in   $H^s(D)$.

\bl\label{liso}
If $D\subset\RR^2$ is a two-dimensional domain with Lipschitz boundary,
then for any $h\in (\tilde H^\frac12(D))'$, the unique solution $w\in H^1_0(D)$
of
\[
-\Delta w= h \hbox{ in } D,
\]
  belongs to $H^{1+\varepsilon}(D)\cap H^1_0(D)$ for any $\varepsilon\in (0,\frac12)$.
\el

\begin{proof}
We fix $\varepsilon\in (0,\frac12)$. Since
$H^{1-\varepsilon}_0(D)=\tilde H^{1-\varepsilon}(D)$ is continuously and densely embedded into
$\tilde H^\frac12(D)$, by duality we obtain
that
$(\tilde H^\frac12(D))'$ is continuously embedded into
$(H^{1-\varepsilon}_0(D))'=H^{-1+\varepsilon}(D)$.
Hence, $w$ can be seen as a solution of the Laplace equation with datum in
$H^{-1+\varepsilon}(D)$.
Owing to Theorem 18.13 and Remark 18.17/2 in \cite{dauge:88}, $w$ belongs to
$H^{1+\varepsilon}(D)$.
\end{proof}

\bl\label{lregpartu}
Let $u\in V$ be the solution of \eqref{weakform}, then we have
\[
\partial_\nu u\in L^2(\Gamma_V).
\]
\el

\begin{proof}
We first observe that, by Theorem 2.8 of \cite{Nicaise92b}, $\partial_\nu u\in
(\tilde H^\frac12(\Gamma_V))'$.
Then,  we may interpret  $u_F\in H^1_0(F)$ as the unique variational
solution of
\be\label{pbonF}
\Delta_{LB} u_F=-g+\partial_\nu u\in (\tilde H^\frac12(F))'.
\ee
By Lemma \ref{liso}, we deduce that
$u_F$ belongs to $H^{1+\varepsilon}(F)\cap H^1_0(F)$ for any $\varepsilon\in (0,\frac12)$.
We now fix $\varepsilon\in (0,\frac12)$ small enough  that
the mapping
\[
H^{\frac{3}{2}+\varepsilon}(\Om)\cap H^1_0(\Om) \to H^{-\frac{1}{2}+\varepsilon}(\Om):
v\to \Delta v,
\]
is an isomorphism (see \cite[Thm 18.13]{dauge:88}).

Now, by applying the trace theorem from \cite{grisvard:74}, there exists $w\in
H^{\frac{3}{2}+\varepsilon}(\Om)$
 such that
 \beq\label{Dirichletbc}
  \gamma w=0\, \hbox{ on } \Gamma_D,
  \\
 \gamma w=u_F\, \hbox{ on } F.
 \label{traceonF}
 \eeq
This implies, again by uniqueness, that $v:=u-w\in H^1_0(\Om)$ is the solution
of
\[
-\Delta v=f+\Delta w\in H^{-\frac{1}{2}+\varepsilon}(\Om).
\]
We therefore deduce that
$v$ belongs to $H^{\frac{3}{2}+\varepsilon}(\Om)$ and, hence, $u$ belongs to
this space  as well.
By a standard trace theorem, we finally conclude that
$\nabla u\in H^{\varepsilon}(\Gamma)^3$.
\end{proof}

Thus, we have the following decomposition of the singular solution $u_F$ on the polygonal face.

\bc\label{creguF}
Let again $u\in V$ be the solution of  \eqref{weakform}. Then, it holds
\be\label{decompuF}
u_F=u_{F,R}+\sum_{S\in {\mathcal S}_F: \omega_{F,S}>\pi} c_S r_S^{\frac{\pi}{\omega_{F,S}}} \sin \left(\frac{\pi \theta_S}{\omega_{F,S}}\right),
\ee
where $u_{F,R}\in H^2(F)$ and $c_S\in \C$.
\ec

\begin{proof}
As $\partial_\nu u$ belongs to $L^2(F)$ by Lemma \ref{lregpartu}, the right-hand
side in \eqref{pbonF} is now in $L^2(F)$
and therefore   Theorem 4.4.3.7 of \cite{grisvard:85a} yields \eqref{decompuF}.
\end{proof}

We  will refer to $u_{F,R}$ as the {\em regular part} of $u_F$, hence the
subscript, as it has the expected regularity from elliptic theory. We will
consequently call  $u_F-u_{F,R}$ the {\em singular part}  of $u_F$.

 For the regularity of the solution in the interior of the domain $\Omega$, we first have the following lifting estimate based on the trace theorem.

 \bl\label{llifting}
Given
$u_{F}\in H^{\frac{3}{2}+\varepsilon}(F)\cap H^1_0(F)$ for some $\varepsilon\in
(0,\frac12)$, there exists a lifting $w\in H^{2+\varepsilon}(\Om)$
 satisfying
\eqref{Dirichletbc}-\eqref{traceonF}.
\el

\begin{proof}
The idea is to use again the trace theorem from \cite{grisvard:74} with
$s=2+\varepsilon$
and the operator
\be\label{eq.trace}
u\to (f_{j, 0}, f_{j,1})_{j=1}^N:=\Big( u_{|\Gamma_j},  ( \partial_{\nu_j} u)_{|\Gamma_j})\Big)_{j=1}^N,
\ee
where $\Gamma_j$ are the faces of $\Om$ and  $\nu_j$ the outward normal vector along $\Gamma_j$.
As $1+\varepsilon$ is not an integer, this trace operator (\ref{eq.trace}) is surjective from $H^{2+\epsilon}(\Omega)$ onto
the subspace of
$\prod_{j=1}^N ( H^{\frac{3}{2}+\varepsilon}(\Gamma_j)\times  H^{\frac{1}{2}+\varepsilon}(\Gamma_j))$
that satisfies the compatibility conditions $(C_1)$ of \cite{grisvard:74}.
If we assume that $\Gamma_1=F$, it is therefore sufficient to show that
there exist $f_{j,1}\in H^{\frac{1}{2}+\varepsilon}(\Gamma_j)$, $j=1,\cdots, N$ such that
\[
(u_F, f_{1,1})\times (0, f_{j,1})_{j=2}^N
\]
satisfies these conditions $(C_1)$.
Since such conditions are quite technical to check, as in \cite{grisvard:74}
we can  reduce to check such conditions in the case where $\Om$ is the trihedral
$x_i>0$, $i=1,2,3$
and $F$ is the face $x_1=0$ (and hence $N=3$ with $\Gamma_2\equiv x_2=0$
and $\Gamma_3\equiv x_3=0$), by means of a localization argument and a linear
change of variables. In such a case, the conditions $(C_1)$
of \cite{grisvard:74} for $(u_F, f_{1,1})\times (0, f_{j,1})_{j=2}^3$
take the form:
\begin{subequations} \label{C1cond}
\beq\label{C1-1}
u_F=0 \hbox{ on } A_{1,2}\cup A_{2, 3},
\eeq
\beq
 f_{1,1}=0 \hbox{ on } A_{1,3},
 \label{C1-2}
\eeq
\beq
 \partial_2 u_F=0 \hbox{ on } A_{1,3},
\label{C1-3}
\eeq
\beq
 \partial_3 u_F= f_{3,1} \hbox{ on } A_{1,3},
 \label{C1-4}
\eeq
\beq
  f_{1,1}=0 \hbox{ on } A_{1,2},
  \label{C1-5}
\eeq
\beq
 \partial_2 u_F=f_{2,1} \hbox{ on } A_{1,2},
 \label{C1-6}
\eeq
\beq
 \partial_3 u_F=0 \hbox{ on } A_{1,2},
 \label{C1-7}
\eeq
\beq
f_{2,1}=0 \hbox{ on } A_{2,3},
\label{C1-8}
\eeq
\beq
f_{3,1}  =0 \hbox{ on } A_{2,3},
\label{C1-9}
\eeq
\end{subequations}
where $A_{i,j}=\bar \Gamma_i\cap \bar \Gamma_j$. The first condition trivially holds as $u_F$ belongs to $H^1_0(F)$,
and similarly \eqref{C1-3} (resp. \eqref{C1-7}) because  $\partial_2 u_F$ (resp. $\partial_3 u_F$) is the tangential derivatives
of $u_F$ on $A_{1,3}$ (resp. $A_{1,2}$).
To satisfy  the second and fourth conditions we simply take  $f_{1,1}=0$. Hence
it remains to verify the conditions
\eqref{C1-6} and \eqref{C1-8} (resp. \eqref{C1-4} and \eqref{C1-9}) that can be
interpreted as constraints on $f_{2,1}$
and $f_{3,1}$, respectively. In other words, we look for $f_{2,1}\in
H^{\frac{1}{2}+\varepsilon}(\Gamma_2)$
(resp. $f_{3,1}\in H^{\frac{1}{2}+\varepsilon}(\Gamma_3)$) satisfying the boundary conditions \eqref{C1-6} and \eqref{C1-8}
(resp. \eqref{C1-4} and \eqref{C1-9}). Such a solution $f_{2,1}$ (and similarly
$f_{3,1}$)  exists by applying Theorem 1.5.1.2 of \cite{grisvard:85a} (valid for
a quarter plane), because the function defined by $ \partial_2 u_F$ on $A_{1,2}$
and 0 on $A_{2,3}$ belongs to $ H^{\varepsilon}(\Gamma_2)$.
\end{proof}

For a vertex $v\in\mathcal S$, let $C_v$ be the infinite polyhedral cone that
coincides with $\Omega$ in the neighborhood of $v$. Let $G_v=C_v\cap S^2(v)$ be
the intersection of $C_v$ and the unit sphere centered at $v$. For an edge
$e\in\mathcal E$, let $\omega_e$ be the interior angle between the two faces of
$\Omega$ that contain $e$. Then, for $v\in\mathcal S$ and for $e\in\mathcal E$,
respectively, we define the following parameters associated to the singularities
in the solution near $v$ and $e$:
\be\label{eqn.ev}
\lambda_{v}:=-\frac{1}{2}+\sqrt{\lambda_{v, 1}+\frac{1}{4}},\qquad\qquad \quad \lambda_e:=\pi/\omega_e,
\ee
where $\lambda_{v,1}$ is the smallest positive eigenvalue of the
Laplace-Beltrami operator on $G_s$ with Dirichlet boundary conditions. We
observe that a vertex $v$ is singular if $\lambda_v<1/2$ and an edge $e$ is
singular if $\lambda_e<1$. For a given macro element $\Lambda_\ell$, we set
$\lambda_v^{(\ell)}=\lambda_v$ if $\Lambda_\ell$ contains one singular vertex
$v$ of $\Omega$ and $\lambda_v^{(\ell)}=\infty$ otherwise. Similarly,  we set
$\lambda_e^{(\ell)}=\lambda_e$ if $\Lambda_\ell$ contains on singular edge $e$
of $\Omega$ and $\lambda_e^{(\ell)}=\infty$ otherwise. Then,  the
following decomposition for the variational solution of (\ref{Ventcel}) holds.

\bt\label{tregu}
Let $u\in V$ be again the solution of \eqref{weakform}. We have:
\be\label{regu}
u=u_R+u_S,
\ee
where $u_{R}\in H^2(\Omega)$ and  $u_S\in H^1(\Om)$ satisfies, for all $\ell\in \{1,\ldots, L\}$,
\beq\label{2.26}
\frac{\partial u_S}{\partial x_j^{(\ell)}}\in V^1_{\beta,\delta}(\Lambda_\ell), j=1,2,\\
\frac{\partial u_S}{\partial x_3^{(\ell)}}\in V^1_{\beta,0}(\Lambda_\ell),
\label{2.27}
\eeq
for any $\beta, \delta\geq 0$ such that
\[
\beta>\frac12-\lambda_v^{(\ell)}, \qquad  \delta>1-\lambda_e^{(\ell)},
\]
\et

Again, the subscripts refer to the fact that $u_R$ has the expected
regularity, and hence it will be called the {\em regular part} of the solution,
while $u_S= u- u_R$ represents the {\em singular part}.

\begin{proof}
The decomposition  \eqref{decompuF} implies that there exists $\varepsilon\in (0,\frac12)$ small enough such that
$u_{F}\in H^{\frac{3}{2}+\varepsilon}(F)\cap H^1_0(F)$.
Hence by Lemma \ref{llifting}, there exists a lifting $w\in H^{2+\varepsilon}(\Om)$ satisfying
\eqref{Dirichletbc}-\eqref{traceonF}.
With this lifting at hands, we consider $u-w$,  which belongs to $H^1_0(\Om)$
and is the weak solution of
\[
-\Delta (u-w)=f+\Delta w.
\]
As $f+\Delta w$ belongs to $L^2(\Om)$, we can apply Theorem 2.10 of
\cite{apel:97c} to $u-w$, which gives the decomposition:
\[
u-w=u_R+u_S,
\]
with $u_R\in H^2(\Omega)$ and $u_S$ satisfying \eqref{2.26}-\eqref{2.27}.
Finally, the result follows by setting $u_R=w+u_r$.
\end{proof}

{\rk Theorem \ref{tregu} shows that for the solution to (\ref{weakform}) with the Ventcel boundary condition, its regularity in $\Omega$, determined by the geometry of the domain, is similar to the regularity of the Poisson equation with the Dirichlet boundary condition. Meanwhile, the trace of the solution $u$ on the face $F$ is the solution of a two-dimensional elliptic problem with the Dirichlet boundary condition.  Corollary \ref{creguF} implies that the regularity of the trace depends on the interior angles of the polygon $F$.}

\section{Finite element approximation\label{sapprox}}

We consider an (anisotropic) triangulation $\caT_h=\{T_i\}_{i=1}^N$ of $\Omega$
as in Section 3 of   \cite{apel:97c} or  in Section 2 of \cite{Apeletal:14},
consisting of tetrahedra
with refinement parameters $\mu_\ell$ and $\nu_\ell$.
We assume the general conditions for a triangulation of the domain
(see e.g.\cite{BrennerScott,Ciarlet}) and  that the number of
tetrahedra $m$ satisfies $N\sim h^{-3}$, where $h$ is the global mesh size.
In addition, we assume that the initial subdomains $\Lambda_\ell$ are resolved
exactly, namely, $\bar\Lambda_\ell=\cup_{i\in L_\ell}\bar T_i$, where $\ell=1,
\cdots, L$
and $L_\ell\subset\{1,\cdots, m\}$ is the index set of the tetrahedra included
in $\bar\Lambda_\ell$.

In each $\Lambda_\ell$, the parameters
$\mu_\ell,
\nu_\ell\in(0, 1]$ determine the anisotropic mesh refinement close to edges and
vertices, respectively as indicated in \eqref{eqn.h} below.
When $\mu_\ell=1$ or $\nu_\ell=1$, there will be no graded refinement
in $\Lambda_\ell$.
We recall the local Cartesian
coordinate system $(x_1^{(\ell)}, x_2^{(\ell)},x_3^{(\ell)})$ in each of the
subdomain $\Lambda_\ell$, which is such that the singular vertex is at the
origin and the singular edge is along the $x_3$-axis, if they exist.  Then, for
each element $T_i\subset\Lambda_\ell$ of the triangularization, we let
$$
r_i:=\inf_{x\in T_i}[(x_1^{(\ell)})^2+(x_2^{(\ell)})^2]^{1/2}, \qquad
R_i:=\inf_{x\in T_i}[(x_1^{(\ell)})^2+(x_2^{(\ell)})^2+(x_3^{(\ell)})^2]^{1/2},
$$
be the distance of $T_i$ to the origin and the $x_3$-axis, respectively.
We then introduce local, anisotropic mesh parameters in $T_i$ as follows:
\begin{eqnarray}\label{eqn.h}
h_i:=\left\{\begin{array}{ll}
h^{1/\mu_\ell}\quad {\rm{if}}\ r_i=0,\\
hr^{1-\mu_\ell}_i \quad {\rm{if}}\ r_i>0,
\end{array}\right.
\qquad H_i:=\left\{\begin{array}{ll}
h^{1/\nu_\ell}\quad {\rm{if}}\ 0\leq R_i\lesssim h^{1/\nu_\ell},\\
hR^{1-\nu_\ell}_i \quad {\rm{if}}\ R_i\gtrsim h^{1/\nu_\ell},
\end{array}\right.
\end{eqnarray}
We  also introduce the actual mesh sizes $\tilde h_{j, i}$, which are
the lengths of the projections of $T_i\subset\Lambda_\ell$ on the
$x^{(\ell)}_j$-axis, $1\leq j\leq 3$. Then, there exists a triangulation
$\mathcal T_h$ satisfying the following conditions:
\begin{enumerate} [label=\arabic*.]
\item If $\mu_\ell<1$, then $\tilde h_{j, i}\sim h_i$, $j=1, 2$, $\tilde
h_{3, i}\lesssim H_i$, and $\tilde h_{3, i}\sim H_i$ if $r_i=0$.
\item The number of tetrahedra in $\Lambda_\ell$ with $r_i=0$ is of order
$h^{-1}$.
\item The number of tetrahedra in $\Lambda_\ell$ such that $0\leq
R_i\lesssim h^{1/\nu_\ell}$ is bounded by $h^{2\mu_\ell/\nu_\ell-2}$, and there
is only one tetrahedral element with $R_i=0$.
\item If $\mu_\ell<1$, then $\mu_\ell\leq \nu_\ell$ for $1\leq \ell\leq L$.
\end{enumerate}
We refer to \cite{apel:97c} or a detailed description of these conditions. It is clear that this triangulation $\mathcal T_h$ induces an exact
triangulation $\mathcal F_h$ of the face $F$, the elements of which are simply given
by $\bar T\cap \bar F$ for  $T\in \caT_h$.

Based on these triangulations, we introduce the approximation space $V_h$ of $V$ as follows:
\[
V_h:=\{u_h\in V: u_{|T}\in \mathbb{P}_1(T),\ \forall T\in \caT_h\},
\]
where $\mathbb{P}_m$, $m\in \ZZ_+$, denotes the space of all polynomials of
degree $\leq m$.
This is clearly a closed subspace of $V$.

Then, the Finite Element approximation of Problem \eqref{weakform} consists of
looking for a solution
$u_h\in V_h$ of
\be\label{weakformh}
a(u_h,v_h)=\int_\Om f  v_h\, dx+\int_{\Gamma_V} g v_h\, d\sigma(x),
\qquad\forall  v_h\in V_h.
\ee
By C\'ea's lemma, we have
\[
\|u-u_h\|_V\lesssim \inf_{v_h\in V_h} \|u-v_h\|_V,
\]
Hence an error estimate will be available if we can built an appropriate approximation $v_h$
of $u$. This is the purpose of the next theorems in this section.

\bt\label{tinterpolation}
Recall the parameters in (\ref{eqn.ev}) and in (\ref{eqn.h}). Assume that for
all $\ell=1,\cdots, L$, we have:
\begin{subequations}
\beq
\mu_\ell <\lambda_e^{(\ell)},
\label{5.2}
\eeq
\beq
\nu_\ell <\lambda_v ^{(\ell)}+\frac12
\label{5.3}
\eeq
\beq
\frac{1}{\nu_\ell} +\frac{1}{\mu_\ell}(\lambda_v ^{(\ell)}-\frac12)>1.
\label{5.4}
\eeq
\end{subequations}
Then, there exists $v_h\in V_h$ such that
\be\label{04/04:1}
  \|u-v_h\|_{1,\Omega}\lesssim h.
\ee
\et

\begin{proof}
The proof of Theorem \ref{tregu} furnishes the splitting
of $u$ as
\[
u=\tilde u+w,
\]
where $\tilde u=u-w$ and $w\in H^{2+\varepsilon}(\Omega)$
with $\varepsilon\in (0,\frac12)$. Hence we define
an interpolant $I_hu$ of $u$
as follows:
\begin{equation} \label{eq.uInterp}
I_hu:=\tilde u_I+D_h  (\tilde u-\tilde u_I)+ L_h w,
\end{equation}
where $D_h$ is the interpolant introduced in  \cite{Apeletal:14},
$\tilde u_I$ is the Lagrange interpolant of $\tilde u$ with respect to the
partition $\{\Lambda_\ell\}$,
while $L_hw$
is the standard Lagrange interpolant of $w$, which consists of piece-wise
polynomials of degree $1$.
Then, using the regularity estimate in Theorem \ref{tregu} and applying Theorem 3.11 of \cite{Apeletal:14}, we have
\be\label{04/04:2}
\|\tilde u-\tilde u_I+D_h  (\tilde u-\tilde u_I)\|_{1,\Omega}\lesssim h.
\ee
On the other hand, as $w$ belongs to $H^{2+\varepsilon}(\Omega)$
and $H^{2+\varepsilon}(\Omega)$ is continuously embedded into $W^{2,p}(\Omega)$
with $p\in (2, \frac{6}{3-2-\varepsilon})$,
by the estimate (5.6) of  \cite{apel:97c} for a fixed $p\in (2, \frac{6}{3-2-\varepsilon})$, we deduce  that
\[
\|w-L_hw\|_{1,\Omega}\lesssim h \|w\|_{W^{2,p}(\Omega)}\lesssim  h \|w\|_{2+\varepsilon,\Omega}.
\]
This estimate and \eqref{04/04:2} prove the estimate \eqref{04/04:1}.
\end{proof}

We observe that $I_h u=L_h w=L_hu_F$ on the face $F$, since $\tilde u_I$ and  $D_h  (\tilde u-\tilde u_I)$ vanish on $F$. We next state and prove
an error estimate for the Finite Element
approximation on the face $F$.

\bt\label{tinterpolation2D} For a macro element $\Lambda_\ell$ such that
$\bar\Lambda_\ell\cap F\neq \emptyset$, let $\omega_{F, v, \ell}$ be the
interior angle of $F$ associated with the vertex $v\in\mathcal V$.
Assume that the conditions
\begin{subequations}
\beq\label{5.32D}
\nu_\ell <\frac{\pi}{\omega_{F,v,\ell}},
\eeq
\beq
\frac{1}{\nu_\ell} +\frac{1}{\mu_\ell}(\frac{\pi}{\omega_{F,v,\ell}}-1)>1,
\label{5.42D}
\eeq
\end{subequations}
are satisfied. Then, it holds:
\be\label{05/04:1}
  \|u_F-L_h u_F\|_{1,F}\lesssim h.
\ee
\et

\begin{proof}
We will prove that for all $\ell=1,\cdots, L$, we have
\be\label{05/04:2}
  |u_F-L_h u_F|_{1,F\cap \bar\Lambda_\ell}\lesssim h.
\ee
Hence, summing on $\ell$, we find that
\[
  |u_F-L_h u_F|_{1,F}\lesssim h,
\]
and the conclusion of the theorem follows from Poincar\'e's inequality.

\noindent  To prove \eqref{05/04:2}, we distinguish different cases: \begin{enumerate}[label=\arabic*.]
\item $\bar F\cap \bar \Lambda_\ell$ contains no singular vertex or singular edge: In this case,  $u_F$ belongs to $H^2(F\cap
\bar\Lambda_\ell)$
and the mesh on $F\cap  \bar\Lambda_\ell$ is quasi-uniform. Thus,  the estimate \eqref{05/04:2} is standard.

\item  $\bar F\cap \bar \Lambda_\ell$ contains a singular vertex $v$ but no singular edge: Thus, $u_F$ belongs to
\[
V^2_\gamma(F\cap  \bar\Lambda_\ell)=\{v\in L^2_{\rm loc}(F\cap  \bar\Lambda_\ell):
R^{\gamma+|\beta|-2} D^\beta v\in L^2(F\cap  \bar\Lambda_\ell), \forall |\beta|\leq 2\},
\]
with $\gamma>1-\frac{\pi}{\omega_{F,v,\ell}}$,
and the estimate \eqref{05/04:2} is also standard, since the triangulation in
$F\cap \bar\Lambda_\ell$ is isotropic (see for instance
\cite{raugel:78}, \cite[\S 8.4]{grisvard:85a}).

\item   $\bar F\cap \bar \Lambda_\ell$ contains a singular edge: Then, the mesh on $F\cap \bar\Lambda_\ell$ is anisotropic.  There are two possible situations: (S1) $\bar F\cap \bar \Lambda_\ell$ contains no singular vertex; and (S2) $\bar F\cap \bar \Lambda_\ell$ also contains a singular vertex $v$. Due to Corollary \ref{creguF},  for (S1), $u_F$ belongs to
$H^2(F\cap \bar\Lambda_\ell)$, while for (S2), $u_F$ belongs to $V^2_\gamma(F\cap \bar\Lambda_\ell)$.   Now for any triangle  $T_i$  in $F\cap \bar\Lambda_\ell$, we will prove that
\be\label{06/04:1}
|u_F-L_hu_F|_{1,T_i}\lesssim h |u_F|_{2,\gamma, T_i},
\ee
with $\gamma=0$ for (S1)  and $\gamma>1-\frac{\pi}{\omega_{F,v,\ell}}$ for
(S2), where
\[
|u|_{2,\gamma, T}^2=\sum_{|\alpha|=2}\int_T R^{2\gamma} |D^\alpha u|^2\,dx.
\]
If this estimate is valid, then summing on $T_i$, we get \eqref{05/04:2}.
\end{enumerate}

\noindent To prove \eqref{06/04:1}, we distinguish two cases.
\begin{enumerate}[label=\roman*.]
\item If $T_i$ is far from the singular corner, then we know that $u_F$ belongs
to $H^2(T_i)$,
and, by using Estimate \eqref{4.92D} below, we have:
\beq
\label{08/04:1}
|u_F-L_hu_F|_{1,T_i}&\lesssim& h_i |\partial_1 u_F|_{1,T_i}+ H_i |\partial_3 u_F|_{1,T_i}
\\
&\lesssim& H_i |u_F|_{2,T_i}.
\nonumber
\eeq
If $\Lambda_\ell$ is of Type 3, then   $u_F$ belongs to $H^2(F\cap\bar\Lambda_\ell)$, but
as $\nu_\ell=1$, by the assumptions on the mesh
we have $H_i\lesssim h$, hence the estimate (\ref{08/04:1}) directly yields
\eqref{06/04:1}.
If $\Lambda_\ell$ is of Type 4,
 we again distinguish two cases:
\begin{enumerate}[label=\alph*)]
 \item
 If $R_i\gtrsim h^\frac{1}{\nu_\ell}$, then $H_i
\sim  h R_i^{1-\nu_\ell}$, and therefore,
\[
|u_F-L_hu_F|_{1,T_i}\lesssim   h R_i^{1-\nu_\ell} |u_F|_{2,T_i}.
\]
This yields \eqref{06/04:1} by our assumption \eqref{5.32D}.

\item If $0<R_i\lesssim h^\frac{1}{\nu}$, then as $H_i\sim  h^\frac{1}{\nu}$,
the estimate \eqref{08/04:1}
becomes:
\[
|u_F-L_hu_F|_{1,T_i}\lesssim  h^\frac{1}{\nu} |u_F|_{2,T_i}.
\]
But from Lemma \ref{sergetriang} below,
we know that $R_i\gtrsim h^\frac{1}{\mu}$ and, therefore,
\[
|u_F-L_hu_F|_{1,T_i}\lesssim   h^{\frac{1}{\nu}-\frac{\gamma}{\mu}} R_i^{\gamma} |u_F|_{2,T_i}.
\]
This yields \eqref{06/04:1} by our assumption \eqref{5.42D}.
\end{enumerate}

\item If $T_i$ is near a singular corner (i.e., $R_i=0$), then applying Lemma
\ref{raugelanisotropic}
we have:
\[
|u_F-L_hu_F|_{1,T_i}\lesssim   h^{\frac{1}{\nu}-\frac{\gamma}{\mu}}   |u_F|_{2,\gamma,T_i}.
\]
Again we get  \eqref{06/04:1} owing to our assumption \eqref{5.42D}.
\end{enumerate}
The proof is now complete.
\end{proof}

Theorems \ref{tinterpolation} and \ref{tinterpolation2D} directly lead to the
following {\em a priori} global interpolation estimate on $u$
and error estimate on the Finite Element solution $u_h$.

\bc\label{cinterpolation}
Assume that for all $\ell=1,\cdots, L$, \eqref{5.2}, \eqref{5.3}, \eqref{5.4},   \eqref{5.32D} and \eqref{5.42D} hold.
Then there exists $v_h\in V_h$ such that
\[
  \|u-v_h\|_V\lesssim h.
\]
\ec

\bc\label{cerrorestimate}
Under the assumption of Corollary \ref{cinterpolation},
if $u\in V$ is the solution of \eqref{weakform} and $u_h\in V_h$ the solution of \eqref{weakformh},
then
\be
  \|u-u_h\|_V\lesssim h.\nonumber
\ee
\ec

\section{Anisotropic error estimates in two dimension\label{sinterror2D}}

To complete the  proof of Theorem \ref{tinterpolation2D}  we need some
interpolation estimates in two space dimensions. In this section, $T_i$ will be
a triangle in the triangulation $\mathcal F_h$ of the face $F$, which is induced
by the triangulation $\mathcal T_h$ of $\Omega$.
We will need  the two-dimensional version of Theorem
4.10
of \cite{apel:97c}, given below.

\bt
Assume that $\Lambda_\ell$ is of Type 3 or 4. Suppose that $\bar F\cap \bar\Lambda$ contains the singular edge. Recall the local Cartesian coordinate system $(x_1, x_2, x_3)$ for $\Lambda_\ell$, for which the singular edge is on the $x_3$-axis. Let $F\cap \bar\Lambda$ be in the plane given by $x_2=0$. Let $T_i\subset  F\cap\bar\Lambda_\ell$ be a triangle in the triangulation $\mathcal F_h$.  Then, for $v\in H^2(T_i)$,  we have
\be\label{4.92D}
|v-L_hv|_{1,T_i}\lesssim h_i |\partial_{1} v|_{1,T_i}+ H_i |\partial_3 v|_{1,T_i},
\ee
where  $h_i$ and $H_i$ are defined in (\ref{eqn.h}).
\et

\begin{proof} Let $\tilde h_{1, i}$ and $\tilde h_{3, i}$ be the lengths of the projections of $T_i$ on the $x_1$- and $x_3$-axis, respectively.
We distinguish between the case $\tilde h_{3,i}\lesssim h_i$ or not.
\begin{enumerate}[label=\arabic*.]
\item If $\tilde h_{3,i}\lesssim h_i$, then diam $T_i\lesssim h_i$ (see the
proof of Theorem 4.10
of \cite{apel:97c}) and owing to Theorem 2 of \cite{apel:92}, we have
\[
|v-L_hv|_{1,T_i}\lesssim h_i |v|_{2,T_i},
\]
and \eqref{4.92D} holds since $h_i\lesssim H_i$.

\item If $\tilde h_{3,i}\gtrsim h_i$, then Theorem 1 of \cite{apel:92} on the
reference element and Lemma  4.8
of \cite{apel:97c} yield
\[
|v-L_hv|_{1,T_i}\lesssim \tilde h_{i, 1}|\partial_{1} v|_{1,T_i}+ \tilde h_{3,i} |\partial_3 v|_{1,T_i}.
\]
This estimate implies \eqref{4.92D}, because $\tilde h_{i,1}\lesssim h_i$
and $\tilde h_{3,i}\lesssim H_i$ (see assumption (B1) in \cite{apel:92},
recalling that $\mu_\ell<1$ if a macro element is of Type 3 or 4).
\end{enumerate}
\end{proof}

We continue with an anisotropic error estimate in weighted Sobolev spaces
(compare with  Theorem 1 of \cite{apel:92} for   two-dimensional  triangles in
standard Sobolev spaces and Theorem 4.5
of \cite{apel:97c} for three-dimensional tetrahedra in weighted Sobolev spaces).

\bt
\label{traugelanisotropic}
Let $\hat T$ be the standard reference element of vertices $(0,0), (1,0)$
and $(0,1)$. Denote by $\hat R$ the distance to $(0,0)$.
Let $0\leq \gamma<1$.
Then for all   $u\in V^2_\gamma(\hat T)$,
and $i=1$ or 2,
we have:
\be\label{estraugelanisotropic}
\|\partial_i(u-Lu)\|_{0, \hat T}\lesssim      \|\hat R^\gamma \nabla \partial_i u\|_{0, \hat T},
\ee
where $Lu$ is the Lagrange interpolant of $u$.
\et

\begin{proof}
We first remark that  Lemma 8.4.1.2 of \cite{grisvard:85a}
shows that $V^2_\gamma(\hat T)$ is continuously embedded into
$C(\overline{\hat T})$, hence the   Lagrange interpolant $Lu$ of $u$ is well-defined.
We define the space:
\[
H^{1,\gamma}(\hat T):=\{v\in L^2(\hat T): \hat R^\gamma \nabla v\in L^2(\hat T)^2\},
\]
which is an Hilbert space equipped with its natural norm $\|\cdot\|_{1,\gamma}$.
We will also use the semi-norm:
\[
|v|_{1,\gamma}=\|\hat R^\gamma \nabla v\|_{\hat T}, \quad \forall v\in H^{1,\gamma}(\hat T).
\]
Then by the proof of
Lemma 8.4.1.2 of \cite{grisvard:85a}, we know
that
$H^{1,\gamma}(\hat T)$ is embedded into $W^{1,p}(\hat T)$
for all $1<p<\frac{2}{1+\gamma}$,
and hence
compactly embedded into $L^2(\hat T)$.
The first embedding  and a trace theorem also guarantee that any $v\in
H^{1,\gamma}(\hat T)$ satisfies
\be\label{trace}
v\in L^1(\hat e),\qquad \|v\|_{L^1(\hat e)} \lesssim
\|v\|_{1,\gamma},
\ee
for any edge $\hat  e$ of $\hat T$.
The second embedding implies that
\be\label{estnormseminorm}
\|v\|_{1,\gamma}\lesssim
|v|_{1,\gamma},
\ee
for all $ v\in H^{1,\gamma}(\hat T)$ such that
$\int_{\hat T} v\,dx=0$.

Now we follow the arguments of Lemma 3 and Theorem 1 of \cite{apel:92}.
We will first prove the estimate for $\partial_1$. We observe that \eqref{trace}
implies that the functional
\[
F(v)=\int_{\hat e_1} v(x)\,d\sigma,
\]
where $\hat e_1$ is the edge of $\hat T$ parallel to the $\hat x_1$ axis, is
well defined  and continuous on $H^{1,\gamma}(\hat T)$:
\be\label{contnormF}
|F(v)|\lesssim  \|v\|_{1,\gamma},  \forall v\in H^{1,\gamma}(\hat T).
\ee

Next, we note note that
\[
F(\partial_1 (u-Lu))=0,
\]
if $u\in V^2_\gamma(\hat T)$.
We then define the polynomial $q$ of degree 1 by
\[
q(\hat x_1, \hat x_2)=c \hat x_1,
\]
where
\[
c=2\int_{\hat T} (\partial_1 u)(\hat x)\,d\hat x.
\]
With this choice, we see that
\[
\int_{\hat T} (\partial_1 (u-q))(\hat x)\,d\hat x=0,
\]
and therefore by \eqref{estnormseminorm} we obtain:
\be\label{08/04:2}
\|\partial_1 (u-q)\|_{1,\gamma}\lesssim
|\partial_1 (u-q)|_{1,\gamma}=|\partial_1 u|_{1,\gamma}.
\ee

As $q-Lu$ is linear, $\partial_1 (q-Lu)$ is constant, and we can write
\[
\|\partial_1 (q-Lu)\|_{1,\gamma}\lesssim |F(\partial_1 (q-Lu))|=|F(\partial_1
(q-u))|.
\]
By \eqref{contnormF}, we deduce that
\[
\|\partial_1 (q-Lu)\|_{1,\gamma}\lesssim  \|\partial_1 (q-u)\|_{1,\gamma}.
\]

 This estimate and  the triangle inequality   imply that
\[
\|\partial_1(u-Lu)\|_{0, \hat T}\leq \|\partial_1 (u-Lu)\|_{1,\gamma}\leq \|\partial_1 (u-q)\|_{1,\gamma}+\|\partial_1 (q-Lu)\|_{1,\gamma}
\lesssim  \|\partial_1 (q-u)\|_{1,\gamma}
\]
and the conclusion for $\partial_1$ follows from \eqref{08/04:2}.
The estimate for $\partial_2$ follows in an analogous manner.
\end{proof}

Then, we are ready to derive the interpolation error estimate near a singular corner of  $F$.

\bl
\label{raugelanisotropic}
Assume that $\Lambda_\ell$ is of Type 3 or 4.
Let $0\leq \gamma<1$.
If $T_i$ is near a singular corner (i.e., $R_i=0$), then
for any $u_F\in V^2_\gamma(T_i)$,
we have
\[
|u_F-L_hu_F|_{1,T_i}\lesssim   h^{\frac{1}{\nu_\ell}-\frac{\gamma}{\mu_\ell}}   |u_F|_{2,\gamma,T_i}.
\]
\el

\begin{proof}
The result follows by mapping  $T_i$ to $\hat T$ as in Lemma 4.8 of
\cite{apel:97c},
by using the estimate \eqref{estraugelanisotropic}, and then mapping back to $T_i$
by using the properties (3.2) and (3.3) in  \cite{apel:97c}
and the fact that
$\hat R\lesssim h^{-\frac{1}{\mu_\ell}} R_i$ (see \cite[p. 538]{apel:97c}).
\end{proof}

\br{\rm
If
$T_i$ is isotropic, the previous Lemma is well known and can be found in
\cite{raugel:78} (see also \cite[\S 8.4]{grisvard:85a}).
}
\er

\bl
\label{sergetriang}
Assume that $\Lambda_\ell$ is of Type  4.
Let  $T_i$ be a triangle belonging to  $\bar\Lambda_\ell\cap F$ such that   $R_i>0$,
then
\[
 R_i\gtrsim h^\frac{1}{\mu_\ell}.
\]
\el

\begin{proof}
Without loss of generality, by a relabeling, we can always assume that $T_0$ is
the triangle that contains the singular vertex $v_\ell$.
Then, it has two edges that contain $v_\ell$, the first one is the edge in the
$x_1$-axis
and is of length  $\sim  h^\frac{1}{\mu_\ell}$,
while the other one is of length $\sim  h^\frac{1}{\mu_\ell}$.
Moreover, as the angle between these two edges is independent of the mesh,
the ball
of center $v_\ell$ and radius $c h^\frac{1}{\mu_\ell}$ intersects only $T_0$ by
choosing $c$ small enough. The estimate follows from the definition of the
distance.
\end{proof}

\section{Numerical examples \label{snumerics}}

In this section, we present some numerical examples to illustrate the theory
presented in the previous sections.

\begin{figure}
\center
\includegraphics[scale=0.25]{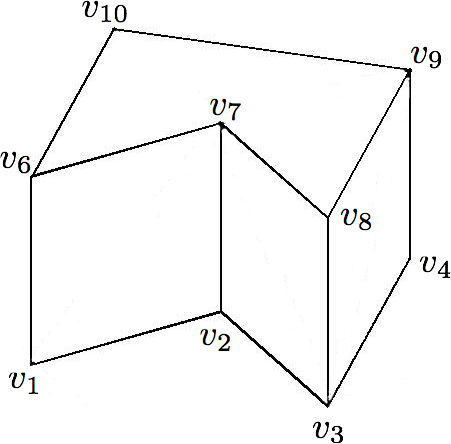}\hspace{2cm}
\includegraphics[scale=0.25]{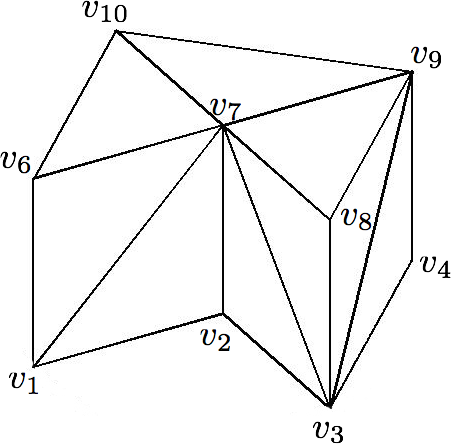}\hspace{0.0cm}
\caption{The computational domain $\Omega$ (left) and the macro elements
(right).}\label{fig.1}
\end{figure}

We will solve the boundary-value problem (\ref{Ventcel}) using the FEM with
linear elements on a polyhedral domain. The domain is given as follows.
We let $\tilde T$ be the triangle with
vertices $(0, 0), (1, 0)$, and $(0.5, 0.5)$, and let the domain be the
prism $\Omega:=\big((0, 1)^2\setminus \tilde T\big) \times(0, 1)$. We refer to
the labeling in Figure \ref{fig.1} in what follows. We will solve
\eqref{Ventcel} in variational form
(\ref{weakform}) with data  $f=1$ and $g=0$. The interior angle between the
two faces that contain the edge $e:=\overline{v_2v_7}$ is $3\pi/2$. Based on the
estimates in (\ref{eqn.ev}) and Theorem \ref{tregu},  $e$ is the  singular edge;
and the solution $u$ admits a decomposition into the singular and regular parts
with regularity determined by $\lambda_e=2/3$. By Theorem \ref{tregu}, the location of the
face $F$, where the Ventcel boundary condition is imposed does not
drastically affect the regularity of the solution.

To verify our theory, we implement two sets of numerical tests regarding different  locations of the special boundary face $F=\Gamma_V$: (I) $F$ is the bottom face of prism $\Omega$, with vertices $v_1, v_2, v_3, v_4$, and $v_5$; (II) $F$ is a face that contains the singular edge with vertices $v_2, v_3, v_7$, and $v_8$.

For both cases, the singular parts $u_S$  of the solution have anisotropic
exponents  and belong to the same weighted space. Moreover, by Corollary
\ref{cerrorestimate}, it is sufficient to  choose the parameters in
(\ref{eqn.h}) corresponding to the singular edge such that  $\mu_\ell<2/3$ and
$\nu_\ell=1$, in order to  achieve the optimal  (first-order) convergence rate.

In Table \ref{tab.5.2new0}, we list the convergence rates of the numerical
solution for the aforementioned model problems with $\nu_\ell=1$, but with
different values of the mesh grading parameter $\mu_\ell$. We let $N$ be the
number of degrees of freedom in the discrete system. Then,  the mesh size
satisfies $h\sim N^{-1/3}$. Since the exact solution is not known, the
convergence rate is computed using the numerical solutions for
successive mesh refinements, $u_{2h}$ $u_{h}$, and $u_{h/2}$,  as
\begin{eqnarray}\label{eqn.rate1}
{\rm{the \ convergence\ rate}}=\log_2(\frac{\|u_{h}-u_{2h}\|_{V}}{\|u_{h/2}-u_{h}\|_{V}}),
\end{eqnarray}
where $u_{2h}$ and $u_{h/2}$ are the finite element solutions with mesh
parameters $2h$ and $h/2$, respectively. Therefore, as $h$ decreases, the
asymptotic convergence rate in (\ref{eqn.rate1}) is a reasonable indicator of
the actual convergence rate for the Finite Element solution.

\begin{table}
\begin{center}
\begin{tabular}{|l|l|||l|}       \hline
\emph{$h\backslash \mu_\ell$} & $0.58$\hspace{0.55cm} $
0.76$\hspace{0.55cm} $1.00$&$0.58$\hspace{0.55cm} $
0.76$\hspace{0.55cm} $1.00$\\ \hline
\hspace{0.1cm}$2^{-3}$ &0.834 \hspace{.26cm} 0.843
\hspace{.25cm} 0.825 \hspace{.26cm}   & 0.821 \hspace{.26cm} 0.833
\hspace{.25cm} 0.825 \hspace{.26cm}\\\hline
\hspace{0.1cm}$2^{-4}$ & 0.938 \hspace{.26cm} 0.930 \hspace{.26cm} 0.890
\hspace{.26cm} & 0.936 \hspace{.26cm} 0.896
\hspace{.26cm} 0.889 \hspace{.26cm}\\ \hline
\hspace{0.1cm}$2^{-5}$ & 0.977 \hspace{.26cm} 0.960 \hspace{.26cm} 0.894
\hspace{.26cm}  & 0.977 \hspace{.26cm} 0.899
\hspace{.26cm} 0.890 \hspace{.26cm} \\ \hline
\hspace{0.1cm}$2^{-6}$ & 0.991 \hspace{.26cm} 0.968 \hspace{.26cm} 0.871
\hspace{.26cm}  & 0.990 \hspace{.26cm} 0.876
\hspace{.26cm} 0.866 \hspace{.26cm}  \\ \hline
\hspace{0.1cm}$2^{-7}$ & 0.995 \hspace{.26cm} 0.968 \hspace{.26cm} 0.837
\hspace{.26cm}  & 1.000 \hspace{.26cm} 0.842
\hspace{.26cm} 0.831 \hspace{.26cm}  \\ \hline
\end{tabular}
\end{center}
\caption{Convergence rates for different values of $\mu_\ell$: (I) $F$ is the bottom face (left); (II) $F$ is a side face containing the singular edge (right).}
\label{tab.5.2new0}
\end{table}

It is clear from the table that for both cases, the first-order convergence rate
is obtained for $\mu_\ell=0.58<2/3$, while we lose the optimal convergence rate
if $\mu_\ell=0.76, 1.00$, both larger than the critical value $2/3$.  When
$\mu_\ell=0.76$, that is, $2/3<\mu_\ell<1$, this choice still leads to an
anisotropic mesh graded toward the singular edge, but the grading is
insufficient to resolve the singularity in the solution, and hence does not give
rise to the predicted first-order convergence rate. These results are in strong
agreement with the theoretical results in Sections \ref{sapprox} and
\ref{sinterror2D}.


\def\cprime{$'$} \def\cprime{$'$} \def\cprime{$'$} \def\cprime{$'$}

\end{document}